\documentclass[12pt]{amsart}
\usepackage{amscd}
\usepackage{latexsym}
\usepackage{amssymb}

\newtheorem{theorem}{Theorem}[section]
\newtheorem{proposition}[theorem]{Proposition}

\newtheorem{lemma}[theorem]{Lemma}

\newtheorem{definition}[theorem]{Definition}

\newtheorem{corollary}[theorem]{Corollary}
\newtheorem{remark}[theorem]{Remark}

\newtheorem{example}[theorem]{Example}

\begin{document}

%
%

\title{The classification of thick representations of simple Lie groups}


\author{Kazunori NAKAMOTO and Yasuhiro OMODA}
\address{Center for Medical Education and Sciences, Faculty of Medicine, 
University of Yamanashi}
\email{nakamoto@yamanashi.ac.jp}
\address{Natural Sciences Division, National Institute of Technology, Akashi College}
\email{omoda@akashi.ac.jp}
\thanks{The first author was partially supported by 
JSPS KAKENHI Grant Number JP23540044, JP15K04814, JP20K03509.}





%
%

\begin{abstract}
We characterize finite-dimensional thick representations over ${\Bbb C}$ of connected complex semi-simple Lie groups by irreducible representations which are weight multiplicity-free and whose weight posets are totally ordered sets. Moreover, using this characterization, we give the classification of thick representations  over ${\Bbb C}$ of connected complex simple Lie groups. 
\end{abstract}

\maketitle

\section{Introduction} 
In our previous paper \cite{Special1}, we have introduced 
$m$-thickness and thickness of group representations. 
Let $\rho : G \to {\rm GL}(V)$ be a finite-dimensional representation of a group $G$.  
If for any subspaces $V_1$ and $V_2$
 of $V$ with $\dim V_1 = m$ and  $\dim V_2 = \dim V - m$  there exists $g \in G$ such that $(\rho(g) V_1) \oplus V_2 =V$, we say that 
a representation $\rho : G \to 
{\rm GL}(V)$ is {\it $m$-thick}. We also say that a representation 
$\rho : G \to {\rm GL}(V)$ is 
{\it thick} if $\rho$ is $m$-thick for each $0 < m < \dim V$ (Definition~\ref{def:thick}). 
In \cite[Proposition~2.7]{Special1}, we proved that 
$1$-thickness is equivalent to irreducibility (Proposition~\ref{prop:imp}).   
Hence $m$-thickness is a natural generalization of 
irreducibility of group representations. 

Let $G$ be a  connected semi-simple Lie group over ${\mathbb C}$, $B$ a Borel
subgroup of $G$, $T$ a maximal torus which is contained in $B$.
Denote their Lie algebras by $\mathfrak{g}, \mathfrak{b}$ and $\mathfrak{t}$,  respectively.
Let $\rho : G \to {\rm GL}(V)$ 
be a finite-dimensional irreducible representation of $G$ over $\mathbb C$.
We denote the set of $\mathfrak{t}$-weights in $V$ by $W(V)$.
Choosing a set of simple roots for $(\mathfrak{g}, \mathfrak{t})$, we can regard $W(V)$ as a partially ordered set (poset) with respect to the usual root order. 
This poset $W(V)$ is called  
the {\it weight poset}.
A representation 
$\rho : G \to {\rm GL}(V)$ is said to be {\it weight multiplicity-free} 
if the weight spaces in $V$ are all one-dimensional.  
We give the following characterization of thickness. 

\begin{theorem}[{Theorem \ref{prop:thick}}] 
An irreducible representation $\rho : G \to {\rm GL}(V)$ 
of a connected semi-simple Lie group $G$ is thick if and only if it is weight multiplicity-free and its 
weight poset is a totally ordered set. 
\end{theorem}

Using this characterization, we can classify the complex thick  
representations of connected semi-simple Lie groups.   
\begin{theorem}[{Theorems~\ref{th:classification} and~\ref{th:classification-semisimple}}] 
If a representation of a connected semi-simple Lie group is thick, then it is geometrically equivalent to one of the following list:
$$e,\ {\rm SL}_{n} (n\ge 2),\  S^{m}{\rm SL}_{2} (m\ge 2),\  {\rm SO}_{2n+1} (n\ge 2),\  {\rm Sp}_{2n} (n\ge 2),\  {\rm G}_{2}. $$
\end{theorem}

\medskip 

Here the irreducible representation of a connected simple Lie group $G$ of the highest weight $\omega_{1}$, where $\omega_{1}$ is the first fundamental weight, is denoted by $G$. 
Similarly, $S^{m}G$ stands for the $m$-th symmetric power of $G$. 
Let $e$ denote the trivial $1$-dimensional representation for any group $G$. 
For the definition of geometric equivalence, see Definition \ref{definition:gm}.

\bigskip

We denote by $\omega_i$ the $i$-th fundamental weight for a connected simple Lie group $G$. 
In \S3, all Lie groups are assumed to be over ${\Bbb C}$ and all representations are finite-dimensional over ${\Bbb C}$.  

The authors would like to express their gratitude to the referee for a detailed review and kind suggestions. 

\section{preliminaries}
A {\it representation} of a group $G$ on a vector
space $V$ is a homomorphism $\rho : G \to {\rm GL}(V)$. 
Then such a map $\rho$ gives $V$ the structure of a $G$-module. 
We sometimes call $V$ itself a representation of $G$ and write $g v$ for $\rho (g)(v)$.  
We recall several definitions and results in our previous
paper \cite{Special1}. 
\begin{definition}[{\cite[Definition~2.1]{Special1}}]\label{def:thick}\rm
Let $G$ be a group. Let $V$ be a finite-dimensional vector space over 
a field $k$. We say that a representation $\rho : G \to {\rm GL}(V)$ is {\it $m$-thick}  
if for any subspaces $V_1$ and $V_2$ of $V$ with $\dim V_1 = m$ and  $\dim V_2 = \dim V - m$,   
there exists $g \in G$ such that $(\rho(g) V_1) \oplus V_2 =V$. 
We also say that a representation $\rho : G \to {\rm GL}(V)$ is 
{\it thick} if $\rho$ is $m$-thick for each $0 < m < \dim V$.   
\end{definition}

\begin{definition}[{\cite[Definition~2.3]{Special1}}]\rm
Let $G$ be a group. Let $V$ be a finite-dimensional vector space over 
a field $k$. We say that a representation $\rho : G \to {\rm GL}(V)$ is 
{\it $m$-dense} if the induced representation 
$\wedge^m \rho : G \to {\rm GL}(\bigwedge^m V)$ is irreducible. 
We also say that a representation $\rho : G \to {\rm GL}(V)$ is 
{\it dense} if $\rho$ is $m$-dense for each $0 < m < \dim V$.   
\end{definition}

We show several examples. See \cite{Special1} for details. 

\begin{example}[{{\it cf. }\cite[Proposition~6.5]{Special1}}]\label{example:GL}\rm  
Let $V$ be the standard representation of ${\rm SL}_n $ and $V^{\ast}$ the dual
representation of $V$. Then $V$ and $V^{\ast}$ are dense. 
\end{example}

\begin{example}[{\cite[Proposition~6.10]{Special1}}]\label{example:SOeven}\rm  
The standard representation of ${\rm SO}_{2n}$ is $m$-dense for each 
$0 < m < 2n$ with $m\neq n$, but not $n$-thick.  
\end{example}

\begin{example}[{\cite[Proposition~6.11]{Special1}}]\label{example:SO}\rm  
The standard representation of ${\rm SO}_{2n+1}$ is dense.  
\end{example}

\begin{example}[{\cite[Proposition~6.18]{Special1}}]\label{example:SP}\rm  
The standard representation of ${\rm Sp}_{2n}$ is thick, but not
 $m$-dense for each $1<m<2n-1$.
\end{example}

Let $V$ be a finite-dimensional representation of a group $G$. 
For positive integers $i$ and $j$ with $i+ j= \dim V$, let us consider the $G$-equivariant perfect 
pairing $\bigwedge^i V \otimes \bigwedge^j V \stackrel{\wedge}{\longrightarrow} \bigwedge^{\dim V} V \cong k$. 
For a $G$-invariant subspace $W$ of $\bigwedge^i V$, put $W^{\perp} := \{ y \in \bigwedge^j V 
\mid  x\wedge y = 0 \mbox{ for any } x \in W \}$. 
Then $W^{\perp}$ is also $G$-invariant. 
In particular, $\bigwedge^i V$ is irreducible if and only if so is $\bigwedge^j V$.

\begin{proposition}[{\cite[Proposition~2.6]{Special1}}]
Let $V$ be an $n$-dimensional representation of a group $G$. 
For each $0 < m < n$, $V$ is $m$-thick $($resp. $m$-dense$)$ if and only if 
$V$ is $(n-m)$-thick $($resp. $(n-m)$-dense$)$. 
\end{proposition}

\begin{proposition}[{\cite[Proposition~2.7]{Special1}}]\label{prop:imp}
For any finite-dimensional representation $V$ of a group $G$, 
the following implications hold 
for $0 < m < \dim V:$   
\begin{eqnarray*} 
\begin{array}{ccccc} 
\mbox{$m$-dense} & \Longrightarrow & \mbox{$m$-thick}  & &  \\ 
 & & \Downarrow & & \\ 
\mbox{$1$-dense} & \Longleftrightarrow & \mbox{$1$-thick} & \Longleftrightarrow & \mbox{irreducible}. 
\end{array}  
\end{eqnarray*}
\end{proposition}

\begin{corollary}[{\cite[Corollary~2.8]{Special1}}]
For any finite-dimensional representation of a group $G$, 
the following implications hold$:$
\[
dense \Rightarrow thick \Rightarrow irreducible. 
\]
\end{corollary}

\begin{corollary}[{\cite[Corollary~2.9]{Special1}}]\label{cor:dim<=3}
For any representation $V$ of a group $G$ with $\dim V\le 3$, the following 
implications hold$:$
\[
dense \Leftrightarrow thick \Leftrightarrow irreducible. 
\]
\end{corollary}

\begin{definition}[{\cite[Definition~2.10]{Special1}}]\label{def:vector}\rm
Let $V$ be an $n$-dimensional vector space over a field $k$. 
For a $d$-dimensional subspace $V'$ of $V$ with $0 < d < n$, we can consider a point  
$[ \bigwedge^d V' ]$ in the projective space ${\mathbb P}(\bigwedge^d V)$. 
In the sequel, we identify $[ \bigwedge^d V' ]$ with   
a non-zero vector $\bigwedge^d V' \in \bigwedge^d V$ (which is determined 
by $[ \bigwedge^d V' ]$ up to scalar) for simplicity.  
For a vector subspace $W \subset \bigwedge^d V$, we say that $W$ is 
{\it realizable} if $W$ contains a non-zero vector $\bigwedge^d V'$ 
obtained by a $d$-dimensional subspace $V'$ of $V$. 
\end{definition}

We have the following criterion of thickness.

\begin{proposition}[{\cite[Proposition~2.11]{Special1}}]\label{prop:condofm-thick}
Let $V$ be an $n$-dimensional representation of a group $G$. 
For $0 < m <n$,      
$V$ is not $m$-thick if and only if there exist 
$G$-invariant realizable subspaces $W_1 \subseteq \bigwedge^m V$ and $W_2 
\subseteq \bigwedge^{n-m} V$ 
such that $W_1^{\perp} = W_2$. 
%
\end{proposition}

\section{The classification of thick representations of simple Lie groups }
Let $G$ be a connected semi-simple Lie group over the complex number field 
$\mathbb C$, $B$ a Borel
subgroup of $G$, $T$ a maximal torus which is contained in $B$, $B^{-}$ 
a Borel subgroup of $G$ opposite to $B$ relative to $T=B\cap B^{-}$.
Denote their Lie algebras by $\mathfrak{g}, \mathfrak{b}, \mathfrak{t}$ and $\mathfrak{b^{-}}$, respectively.
Let $V$ be a finite-dimensional irreducible representation of $G$ over $\mathbb C$.
We will denote the set of $\mathfrak{t}$-weights in $V$ by $W(V)$.
For any weight $\varphi \in W(V)$, let $V_{\varphi}$ be the $\varphi$-weight space in $V$.
Let $\Pi$ be the set of simple roots and ${\Delta}^{+}$ the set of positive roots for ($\mathfrak{g}, \mathfrak{b}$).
We can regard $W(V)$ as a partially ordered set (poset) with respect to the usual root order. More precisely, $\mu > \gamma$ if and only if $\mu -\gamma $ is a nonzero sum of
simple roots with nonnegative coefficients. In particular, if $\mu
-\gamma $ is a simple root, we say that $\mu$ {\it covers} $\gamma$.  
The partially ordered set $W(V)$ is called the {\it weight poset}.  
Following \cite[\S~4.5]{Howe}, 
we say that a representation $V$ of $G$ is {\it weight multiplicity-free} 
({\it WMF}) if the weight spaces in $V$
are all one-dimensional. Howe \cite{Howe} classified the irreducible representations of connected simple Lie groups which are weight multiplicity-free.

\begin{proposition}\label{prop:WMF} 
If a representation $V$ of $G$ is thick, it is weight multiplicity-free.
\end{proposition}

\proof 
Assume that $V$ is not WMF. Then there exists a weight $\varphi \in
W(V)$ such that the dimension of $V_{\varphi}$ is larger than one.
Let $W^{+}(\varphi)$ be the set of all weights strictly larger than
$\varphi$, and $Y^{+}({\varphi})$  the subspace of $V$ which is spanned by all weight
spaces for weights in $W^{+}(\varphi)$. Because the dimension
of $V_{\varphi}$ is larger than one, we can choose two linear independent
$\varphi$-weight vectors $v$ and $w$. 
Let $W^{a,b}_{\varphi}(+)$ be ${\mathbb C} (av+bw)\oplus
Y^{+}({\varphi})$ for $a,b \in {\mathbb C}$.
The subspace $W^{a,b}_{\varphi}(+)$ is 
$B$-invariant. 
Let $n$ be the dimension of $V$, and $d$ the dimension of
$W^{a,b}_{\varphi}(+)$ for $(a, b) \in {\mathbb C}^2\setminus \{ (0, 0) \}$. 
The elements ${\bigwedge}^{d}W^{a,b}_{\varphi}(+)$ 
for $(a, b) \in {\mathbb C}^2\setminus \{ (0, 0) \}$   
are distinct $B$-eigenvectors in ${\bigwedge}^{d}V$ with the same weight.
Let $U_{\varphi}^{a,b}(+)$ be the
irreducible $G$-submodule in ${\bigwedge}^{d}V$ with the highest weight vector
${\bigwedge}^{d}W^{a,b}_{\varphi}(+)$.   
Let $U_{\varphi}(+)$ be the direct sum $U_{\varphi}^{1,0}(+)\oplus U_{\varphi}^{0,1}(+)
\subset {\bigwedge}^{d}V$. Any irreducible $G$-submodule of $U_{\varphi}(+)$ is equal
to $U_{\varphi}^{a,b}(+)$ for some 
$(a, b) \in {\mathbb C}^2\setminus \{ (0, 0) \}$. 
Hence any irreducible $G$-submodule 
of $U_{\varphi}(+)$ is realizable.

Let $Y^{-}(\varphi)$ be the subspace of $V$ which is spanned by all
weight spaces for weights in $W(V)\setminus \{ W^{+}(\varphi), \varphi
\}$.
We take a basis $\{ v,w,u_{1},\ldots ,u_{s} \}$ for $V_{\varphi}$ which
contains $v,w$. 
Let $W_{\varphi}(-)$ be the subspace of $V$ which is spanned by $\{ w,u_{1},\ldots
,u_{s} \}$ and $ Y^{-}(\varphi)$. 
The subspace $W_{\varphi}(-)$ is invariant under the action of the opposite
Borel subgroup $B^{-}$. The equalities $\dim W_{\varphi}(-)=\dim V
-\dim W^{a,b}_{\varphi}(+)=n-d$ hold 
for $(a, b) \in {\mathbb C}^2\setminus \{ (0, 0) \}$. 
Then $\bigwedge^{n-d}W_{\varphi}(-)$ is a $B^{-}$-eigenvector in
$\bigwedge^{n-d}V$. 
Let $U_{\varphi}(-)$ be the
irreducible $G$-submodule with the lowest weight vector
${\bigwedge}^{n-d}W_{\varphi}(-)$ in ${\bigwedge}^{n-d}V$.   
Obviously, $U_{\varphi}(-)$ is realizable.
The irreducibility of $U_{\varphi}(-)$ shows the irreducibility of
$\bigwedge^{d}V /(U_{\varphi}(-))^{\perp}$. Then
$(U_{\varphi}(-))^{\perp}\cap U_{\varphi}(+) \neq \{ 0\}$ because $U_{\varphi}(+)$ is not
irreducible.
Hence $(U_{\varphi}(-))^{\perp}\cap U_{\varphi}(+)$ contains some realizable
$G$-submodule $U_{\varphi}^{a,b}(+)$. 
Therefore $(U_{\varphi}(-))^{\perp}$ is realizable. 
Putting $W_{1}=(U_{\varphi}(-))^{\perp}$ and $W_{2}=U_{\varphi}(-)$,  we see that 
$V$ is not thick by Proposition \ref{prop:condofm-thick}. Hence if $V$ is thick, then it is WMF. 
\qed

\begin{proposition}\label{prop:TOS} 
If a representation $V$ of $G$ is thick, its weight
 poset $W(V)$ is a totally ordered set.
\end{proposition}

\proof 
Let  $V$ be a thick representation of $G$.  
By Proposition \ref{prop:WMF}, $V$ is WMF. 
For any weight $\phi \in W(V)$, 
let $W^{+}(\phi)$ be the set of all weights strictly larger than
$\phi$, and $Y^{+}({\phi})$ the subspace of $V$ which is spanned by
all weight spaces for weights in $W^{+}(\phi)$. 
Note that the irreducible representation $V$ has a highest weight $\omega$ and that 
each weight of $V$ has the form $\omega - \sum_{i=1}^{l} m_i\alpha_i$ ($m_i \in {\mathbb N}$), where 
$\Pi = \{ \alpha_1, \ldots, \alpha_{l} \}$.  

Suppose that the weight poset $W(V)$ is not a totally ordered set. 
There exists an integer $d > 1$ such that $W(V)$ has the $(d-1)$-st highest weight, but not the $d$-th highest weight. 
Let $\varphi$ be the $(d-1)$-st highest weight, and $\psi_1, \psi_2$ maximal weights in $W(V)\setminus (W^{+}(\varphi)\cup \{ \varphi \})$. 
Then the subset 
$W^{+}(\varphi ) \cup \{\varphi\}$ 
is a totally ordered set, $\varphi$
covers $\psi_{1},\psi_{2}$, and $W^{+}({\psi}_{1})=W^{+}({\psi}_{2})=
W^{+}(\varphi ) \cup \{\varphi\}$. 
Because $V$ is WMF, there exists a unique
$\psi_{i}$-weight vector $v_{i}$  up to scalar for each $i=1,2$. 
Let $W_{\psi_{i}}(+)$ be ${\mathbb C} v_{i}\oplus
Y^{+}({\psi_{i}})$.
The subspaces $W_{\psi_{i}}(+)$ are
$B$-invariant for each $i = 1, 2$. 
Let $n$ be the dimension of $V$. Note that $\dim W_{\psi_{i}}(+)  = d$ for $i=1, 2$. 
The elements ${\bigwedge}^{d}W_{\psi_{1}}(+)$ and ${\bigwedge}^{d}W_{\psi_{2}}(+)$ 
are distinct $B$-eigenvectors with distinct weights in ${\bigwedge}^{d}V$.
Let $U_{\psi_{i}}(+)$ be the
irreducible $G$-submodule of ${\bigwedge}^{d}V$ with the highest weight vector
${\bigwedge}^{d}W_{\psi_{i}}(+)$ for each $i=1,2$.
Then $U_{\psi_{1}}(+)$ and $U_{\psi_{2}}(+)$ are realizable and not isomorphic
 to each other as $G$-modules. 
Let $Y^{-}(\psi_{1})$ be the subspace of $V$ which is spanned by all
weight spaces for weights in $W(V)\setminus (W^{+}(\psi_{1})\cup 
\{ \psi_{1} \})$. 
The subspace $Y^{-}(\psi_{1})$ is invariant under the action of the opposite
Borel subgroup $B^{-}$. The equalities $\dim Y^{-}(\psi_{1})=\dim V
-\dim W_{\psi_{1}}(+) =n-d$ hold.
Then $\bigwedge^{n-d}Y^{-}(\psi_{1})$ is a $B^{-}$-eigenvector
in $\bigwedge^{n-d}V$. 
Let $U_{\psi_{1}}(-)$ be the
irreducible $G$-submodule of 
${\bigwedge}^{n-d}V$ with the lowest weight vector
${\bigwedge}^{n-d}Y^{-}(\psi_{1})$.   
Then $U_{\psi_{1}}(-)$ is realizable.
The irreducibility of $U_{\psi_{1}}(-)$ shows the irreducibility of
$\bigwedge^{d}V /(U_{\psi_{1}}(-))^{\perp}$. Then 
$(U_{\psi_{1}}(-))^{\perp}\cap (U_{\psi_{1}}(+)\oplus U_{\psi_{2}}(+))
\neq \{ 0\}$.
Because $U_{\psi_{1}}(+)$ is not isomorphic to $U_{\psi_{2}}(+)$,
$U_{\psi_{1}}(+)\subset (U_{\psi_{1}}(-))^{\perp}$ or
$U_{\psi_{2}}(+)\subset (U_{\psi_{1}}(-))^{\perp}$. 
In particular, $(U_{\psi_{1}}(-))^{\perp}$ is realizable. 
Putting $W_{1}=(U_{\psi_{1}}(-))^{\perp}$ and $W_{2}=U_{\psi_{1}}(-)$, we see that 
$V$ is not thick by Proposition \ref{prop:condofm-thick}. 
This is a contradiction. Hence $W(V)$ is a totally ordered set. 
\qed 

\bigskip 

Let us denote the Grassmann variety which is the set of all $k$-dimensional subspaces 
of a vector space $V$ by ${\rm Grass}(k,V) (\subset {\mathbb P}(\bigwedge^k V))$.  

\begin{lemma}\label{lem:HWV}
Let $V$ be a representation of $G$, and $W$  a $G$-invariant realizable subspace of 
$\bigwedge^{k} V$.
Then there exists $[v] \in {\mathbb P}(W) \cap {\rm Grass}(k,V)$ 
such that $[v]$ is $B$-invariant.
\end{lemma}

\proof 
Let $X$ be ${\mathbb P}(W) \cap {\rm Grass}(k,V)$. 
Because $W$ is realizable, $X$ is not empty. Note that  
$X$ is $G$-invariant and compact.
We take a $G$-orbit $O$ in $X$ whose dimension is minimal. The orbit $O$
is closed and then compact. There is a parabolic subgroup $P$ of $G$ such
that the orbit $O$ is isomorphic to $G/P$. 
Then there is a point $[v] \in O \subset {\mathbb P}(W) \cap
{\rm Grass}(k,V)$ such that $[v]$ is $B$-invariant. 
\qed 

\begin{lemma}\label{lemma:WMF-TOS} 
Assume that an irreducible representation $V$ of $G$ is weight multiplicity-free, its weight poset $W(V)$ is a totally ordered set $\{ \varphi_{1} >\varphi_{2} > \cdots >\varphi_{n} \}$, and $W$ is a $G$-invariant realizable subspace of $\bigwedge^{k} V$. Let $v_{i}$ be a nonzero vector in the $\varphi_i$-weight space $V_{{\varphi}_{i}}$ $(i=1,2,\ldots ,n)$. Then $W$ contains $v_{1}\wedge v_{2}\wedge \cdots \wedge v_{k}$ and $v_{n-(k-1)}\wedge v_{n-(k-2)}\wedge \cdots \wedge v_{n}$.

\end{lemma}

\proof 
Because $V$ is weight multiplicity-free, $\{ v_{1}, \ldots , v_{n} \}$ is a basis of $V$.
By Lemma \ref{lem:HWV}, there exists $[v] \in {\mathbb P}(W) \cap {\rm Grass}(k,V)$ such 
that $v$ is a highest weight vector of an irreducible subrepresentation of $W$ with respect to $B$.
We can put 
$$
 \begin{array}{rl}
 v= & \ \ (p_{1,1}v_{1}+p_{1,2}v_{2}+\cdots +p_{1,n}v_{n})  \\
  &\wedge (p_{2,1}v_{1}+p_{2,2}v_{2}+\cdots +p_{2,n}v_{n}) \\
  & \ \ \ \ \ \ \ \ \   \ \ \ \ \ \ \ \ \  \ \ \ \ \ 
  \vdots \\
  &\wedge (p_{k,1}v_{1}+p_{k,2}v_{2}+\cdots +p_{k,n}v_{n}) \\
\end{array}
$$ 
up to scalar multiplication, 
where $P=(p_{i,j})$ is in reduced row echelon form. Remark that $P$ is uniquely determined.
Let $X_{\alpha}$ be a root vector for a positive root $\alpha \in {\Delta}^{+}$. 
Then $X_{\alpha}v=0$ holds for any $\alpha \in {\Delta}^{+}$. 
If $p_{1,1}=p_{1,2}=\cdots =p_{1,i}=0$ and $p_{1,i+1}=1$ for $i\ge 1$, there is a positive root 
$\alpha \in {\Delta}^{+}$ such that $X_{\alpha}v_{i+1}$ is $c v_{i}$ for a nonzero constant $c$.
Then $X_{\alpha}v$ is not $0$. This is a contradiction.  So $p_{1,1}=1$.
Similarly, we can show that $p_{22}=\dots =p_{kk}=1$.
Because $v$ is a highest weight vector, for any $t \in \mathfrak{t}$ there is a constant $c$ such that  $t v = c v$.  Then by the uniqueness of $P$ we can show that $p_{ij}=0$ for $i=1,\dots ,k$ and $j=k+1,\dots ,n$. 
Then $v=v_{1}\wedge v_{2} \wedge \cdots \wedge v_{k} \in W$.
A similar argument with respect to $B^{-}$ shows that $v_{n-(k-1)} \wedge v_{n-(k-2)}\wedge \cdots \wedge v_{n} \in  W$. 
\qed

\begin{theorem}\label{prop:thick} 
An irreducible representation $V$ of a connected semi-simple Lie group $G$ is thick if and only if it is weight multiplicity-free and its 
weight poset is a totally ordered set. 
\end{theorem}

\proof
The ``only if'' part can be proved by Propositions \ref{prop:WMF} and \ref{prop:TOS}.  
Let us prove the ``if'' part. 
Let us use the notations in Lemma \ref{lemma:WMF-TOS}.
Assume that $W_1 \subseteq \bigwedge^k V$ and $W_2 
\subseteq \bigwedge^{n-k} V$ are $G$-invariant realizable subspaces.  
By Lemma \ref{lemma:WMF-TOS},  
$v_{1}\wedge v_{2}\wedge \cdots \wedge v_{k} \in W_1$ and $v_{k+1} \wedge v_{k+2}\wedge \cdots \wedge v_{n} \in W_2$. Since 
$(v_{1}\wedge v_{2}\wedge \cdots \wedge v_{k}) \wedge (v_{k+1} \wedge v_{k+2}\wedge \cdots \wedge v_{n}) \neq 0$,  $W_1^{\perp} \neq W_2$. 
By Proposition \ref{prop:condofm-thick}, $V$ is thick.
\qed

\bigskip

By \cite[Theorem~4.6.3]{Howe}, we have Howe's classification of irreducible representations of connected simple Lie groups which are weight multiplicity-free. We also refer to Panyushev's paper \cite[Table~1]{panyushev} for the weight posets of weight multiplicity-free representations. Thus, we have 

\begin{theorem}\label{th:listofthick} 
The thick representations of connected simple Lie groups are those on the following list: 
\begin{enumerate}
\item the trivial $1$-dimensional representation for any groups 
\item $A_n \; (n\ge 1)$    
\begin{itemize} 
\item the standard representation $V$ for $A_n \; (n\ge1)$ with highest weight $\omega_1$ 
\item the dual representation $V^{\ast}$ of $V$ for $A_n \; (n\ge1)$ with highest weight $\omega_n$ 
\item the symmetric tensor $S^{m}(V)$ $(m\ge 2)$ of $V$ for $A_1$ with highest weight $m\omega_1$  
\end{itemize}
\item $B_n \; (n\ge 2)$ 
\begin{itemize}
\item  the standard representation $V$ for $B_n \; (n\ge2)$ with highest weight $\omega_1$ 
\item  the spin representation for $B_2$ with highest weight $\omega_2$ 
\end{itemize} 
\item $C_n \; (n\ge 3)$ 
\begin{itemize}
\item  the standard representation $V$ for $C_n \; (n\ge3)$ with highest weight $\omega_1$ 
\end{itemize} 
\item $G_2$ 
\begin{itemize}
\item the $7$-dimensional representation $V$ for $G_2$ with highest weight $\omega_1$.   
\end{itemize} 
\end{enumerate}
\end{theorem}

\proof 
By Theorem~\ref{prop:thick}, it suffices to list up all irreducible representations which are weight multiplicity-free and whose weight posets are totally ordered sets. Using \cite[Theorem~4.6.3]{Howe} and  \cite[Table~1]{panyushev}, we can obtain the list of thick representations of connected simple Lie groups. 
\qed 

\bigskip 

We also have the list of dense representations: 

\begin{theorem}\label{th:listofdense} 
The dense representations of connected simple Lie groups are those on the following list: 
\begin{enumerate}
\item the trivial $1$-dimensional representation for any groups 
\item $A_n \; (n\ge 1)$    
\begin{itemize} 
\item the standard representation $V$ for $A_n \; (n\ge1)$ with highest weight $\omega_1$ 
\item the dual representation $V^{\ast}$ of $V$ for $A_n \; (n\ge1)$ with highest weight $\omega_n$ 
\item the symmetric tensor $S^{2}(V)$ of $V$ for $A_1$ with highest weight $2\omega_1$  
\end{itemize}
\item $B_n \; (n\ge 2)$ 
\begin{itemize}
\item  the standard representation $V$ for $B_n \; (n\ge2)$ with highest weight $\omega_1$.  
\end{itemize} 
\end{enumerate}
\end{theorem} 

\proof
It suffices to verify whether thick representations in the list of Theorems~\ref{th:listofthick} are dense or not.   
It is well-known that the standard representations $V$ of $A_n$ and $B_n$ are dense. We also see that the dual representation $V^{\ast}$ of $V$ for $A_n$ is dense. (For  $A_n$, see  Example~\ref{example:GL} or  \cite[\S15.2]{FH}.  For  $B_n$, see Example~\ref{example:SO} or \cite[Theorem~19.14]{FH}.)  
By Corollary~\ref{cor:dim<=3}, $S^2(V)$ for $A_1$ is dense since $\dim S^2(V)=3$.  

Conversely, let us show that $S^m(V)$ for $A_1$ is not dense if $m \ge 3$. 
Let $\{ \varphi_1 > \varphi_2 \}$ be the weight poset of the standard representation $V$ of $A_1$. 
The weight poset of $S^m(V)$ is $\{ (m-k)\varphi_1+ k \varphi_2 \mid k=0, 1, 2, \ldots, m \}$. 
Thereby, the weight poset of $\bigwedge^2 S^m(V)$ is $\{ (2m-k_1-k_2)\varphi_1 + (k_1 + k_2) \varphi_2 
\mid 0 \le k_1 < k_2 \le m \}$.  
If $m\ge 3$, then $\dim \bigwedge^2 S^m(V)_{(2m-3)\varphi_1+3\varphi_2} = 2$ for the cases $(k_1, k_2)=(0, 3), (1, 2)$.  
Since $\bigwedge^2 S^m(V)$ is not weight multiplicity-free and any irreducible 
representations $S^{m'}(V)$ of $A_1$ are weight  
multiplicity-free, $\bigwedge^2 S^m(V)$ is not irreducible. Hence $S^{m}(V) \ (m\ge 3)$ is not dense.  
It is well-known that the first fundamental representations of $C_n$ and $G_2$ are not dense.  
(For $C_n$, see Example~\ref{example:SP} or \cite[\S17.2]{FH}.   For $G_2$,  
see \cite[\S22.3]{FH}.)  The spin representation for $B_2$ with highest weight $\omega_2$ is not dense since 
it is equivalent to the first fundamental representation for $C_2$ (for $C_2$, see Example~\ref{example:SP} or \cite[\S16.2]{FH}).  
Therefore, we obtain the list of dense representations. 
\qed

\bigskip

According to  \cite[\S6]{BR}, \cite[\S5]{Knop}, and so on, we  introduce the notion of geometric equivalence for 
simplifying  the classification of thick representations.


\begin{definition}[{{\it cf.} \cite[\S6]{BR}, \cite[\S5]{Knop}}]\label{definition:gm}\rm  
For two representations $\rho : G \to {\rm GL}(V)$ and 
$\rho' : G' \to {\rm GL}(V')$, we say that they are {\it geometrically equivalent} if 
there exists a ${\Bbb C}$-linear isomorphism $f:V\to V'$ such that  
$\rho'(G') = f \rho(G) f^{-1}$. 
\end{definition}

We prove the following proposition which was known in \cite[\S5]{Knop}.

\begin{proposition}[{\cite[\S5]{Knop}}]
Let $G$ be a connected semi-simple Lie group over ${\Bbb C}$. 
Let $\rho^{\ast} : G \to {\rm GL}(V^{\ast})$ be the dual representation of a finite-dimensional irreducible  representation $\rho :G \to {\rm GL}(V)$ over ${\Bbb C}$. Then $\rho$ and $\rho^{\ast}$ are geometrically equivalent. 
\end{proposition}

\proof 
Let ${\frak h}$ be a Cartan subalgebra of the Lie algebra ${\frak g}$ of $G$.  
Fix a set of simple roots $\Pi$ of the root system $\Delta$.  
Let $\lambda $ be the highest weight of  $V$ with respect to $\Pi$ 
and $w_{0}$ the longest element of the Weyl group $W$. 
Then $ -w_{0}(\lambda)$ is the highest weight of   $V^{\ast}$ (see \cite[Excercises~10.9 and 21.6]{Humphreys}). 
Let $\phi':{\frak h} \to {\frak h}$ be the isomorphism whose dual  
$\phi'^{\ast} : {\frak h}^{\ast} \to {\frak h}^{\ast}$ is given by $\mu \mapsto -w_{0}(\mu)$. 
There exists a Lie algebra isomorphism $\phi : {\frak g} \to {\frak g}$ 
extending $\phi'$ (see \cite[Theorem~18.4 (b)]{Humphreys}).  
Take a universal cover $\pi : \widetilde{G} \to G$. The dual representation $\widetilde{\rho}^{\ast}$ 
of  $\widetilde{\rho} = \rho \circ \pi : \widetilde{G} \to {\rm GL}(V)$ can be identified with $\rho^{\ast}\circ \pi$.  
By \cite[Chapter III, \S6, Theorem~1]{Bourbaki-Ch1-3}, there exists an automorphism $\psi : \widetilde{G} \to \widetilde{G}$ such that $d\psi = \phi$.  Since $\widetilde{\rho} \circ \psi$ and $\widetilde{\rho}^{\ast}$ have the same highest weight 
$-w_{0}(\lambda)$, there exists an isomorphism $f : V \to V^{\ast}$ such that  
$(\rho^{\ast}\circ \pi) (\widetilde{g}) = \widetilde{\rho}^{\ast} (\widetilde{g}) = f \circ (\widetilde{\rho}\circ \psi)(\widetilde{g})  \circ f^{-1}$ for any $\widetilde{g} \in \widetilde{G}$. 
Hence $\rho^{\ast}(G) = \widetilde{\rho}^{\ast}(\widetilde{G}) = f ((\widetilde{\rho} \circ \psi)(\widetilde{G})) f^{-1} = f \rho(G) f^{-1}$.  
Therefore $\rho$ and $\rho^{\ast}$ are geometrically equivalent. 
 \qed

\begin{remark}\rm 
Assume that two representations $\rho : G \to {\rm GL}(V)$ and 
$\rho' : G' \to {\rm GL}(V')$ are geometrically equivalent.  Then $\rho$ is thick (resp. dense) if and only if so is $\rho'$. 
\end{remark}

According to \cite[\S3.1]{Kac}, we denote the irreducible representation of a connected simple Lie group $G$ with  highest weight $\omega_{1}$ by $G$. 
Similarly, $S^{m}G$ stands for the $m$-th symmetric power of $G$. 
In addition, let $e$ denote the trivial $1$-dimensional representation for any groups $G$. 
Then we have: 


%


\begin{theorem}\label{th:classification} 
If a representation of a connected simple Lie group is thick, 
then it is geometrically equivalent to one of the following list:
$$e,\ {\rm SL}_{n} (n\ge2),\  S^{m}{\rm SL}_{2} (m\ge 2),\  {\rm SO}_{2n+1} (n\ge 2),\  {\rm Sp}_{2n} (n\ge 2),\  {\rm G}_{2}. $$
If a representation of a connected simple Lie group is dense, then it is geometrically equivalent to one of the following list: 
$$e, \ {\rm SL}_n (n\ge 2),\ S^{2}{\rm SL}_{2}, \ {\rm SO}_{2n+1} (n\ge 2).$$
\end{theorem}

\proof 
The last fundamental representation of $B_2$ with highest weight $\omega_2$ is geometrically equivalent to the first fundamental representation of $C_2$ with highest weight $\omega_1$, that is, ${\rm Sp}_4$. 
By Theorems~\ref{th:listofthick} and \ref{th:listofdense}, we have the classification above. 
\qed 

\bigskip 

Theorem~\ref{th:classification} also shows the list of geometric equivalence classes of thick (or dense) representations of connected semi-simple Lie groups. 

\begin{theorem}\label{th:classification-semisimple}
Any thick representation $V$ of a connected semi-simple Lie group $G$ is geometrically equivalent to 
one of the list in Theorem~\ref{th:classification}. In particular, the list of geometric equivalence classes of thick representations (resp. dense representations) of connected semi-simple Lie groups is the same as that of 
thick representations (resp. dense representations) of connected simple Lie groups in 
Theorem~\ref{th:classification}. 
\end{theorem} 

\proof 
Let $\rho : G \to {\rm GL}(V)$ be a thick representation of a connected semi-simple Lie group $G$. 
Take a universal cover $\pi : \widetilde{G} = G_1 \times G_2 \times \cdots \times G_r \to G$, where 
$G_i$ is  a simply-connected simple Lie group for each $i=1, 2, \ldots, r$. We have a thick representation $\widetilde{\rho} = \rho \circ \pi : \widetilde{G} \to {\rm GL}(V)$. 
Since $V$ is an irreducible representation of $\widetilde{G}$, there exist irreducible representations 
$V_i$ of $G_i$ ($1 \le i \le r$) such that $V \cong V_1\otimes V_2 \otimes \cdots \otimes V_r$ as representations of $\widetilde{G}$. 
By Theorem~\ref{prop:thick}, $V$ is WMF as a representation of $\widetilde{G}$ and the weight poset $W_{\widetilde{G}}(V)$ is a totally ordered set. Here, weights in $W_{\widetilde{G}}(V)$ are with respect to a maximal torus $T = T_1\times T_2 \times \cdots \times T_r$ of $\widetilde{G}$, where $T_i$ is a maximal torus of $G_i$.   The order in $W_{\widetilde{G}}(V)$ is defined with respect to a set $\Pi = \Pi_1 \sqcup  \Pi_2 \sqcup \cdots \sqcup \Pi_r$ of simple roots of $\widetilde{G}$, where $\Pi_i$ is a set of simple roots 
of $G_i$.  
Let $W_{G_i}(V_i)$ be the weight poset (with respect to $T_i$ and $\Pi_i$) of the $G_i$-module $V_i$. 
We can write $W_{\widetilde{G}}(V) = \{ \sum_{i=1}^{r} \psi_i \mid \psi_i \in W_{G_i}(V_i) \}$. 

Suppose that there exists $1 \le i < j \le r$ such that $\widetilde{\rho}(G_i) \neq \{ e \}$ and  $\widetilde{\rho}(G_j) \neq \{ e \}$. Then $\sharp W_{G_i}(V_i) \ge 2$ and $\sharp W_{G_j}(V_j) \ge 2$. 
Choose $\phi_1, \phi_2 \in W_{G_i}(V_i)$ and $\varphi_1, \varphi_2 \in W_{G_j}(V_j)$ such that $\phi_1 > \phi_2$ and $\varphi_1 > \varphi_2$. 
Let $\xi = \sum_{k \neq i, j} \psi_k$ be the sum of the highest weights $\psi_k  \in W_{G_k}(V_k)$ 
for $k \neq i, j$.  
For $\eta_1 = \xi + \phi_1 + \varphi_2, \eta_2 = \xi + \phi_2 + \varphi_1 \in W_{\widetilde{G}}(V)$, 
neither $\eta_1 > \eta_2$ nor $\eta_1 < \eta_2$ holds. This implies that 
$W_{\widetilde{G}}(V)$ is not totally ordered, which is a contradiction.   
Hence, any $G_k$ satisfy $\widetilde{\rho}(G_k) = \{ e \}$ except some $G_i$.   
Since $V_k = {\mathbb C}$ except for $k=i$, the representation $V$ of $\widetilde{G}$ is geometrically equivalent to the representation $V_i$ of $G_i$. 
In particular, the representation $V$ of $G$ is geometrically equivalent to a thick representation $V_i$ of 
a connected simple Lie group $G_i$. 
Therefore, Theorem~\ref{th:classification} also shows the lists of geometric equivalence classes of thick and dense representations of connected semi-simple Lie groups. 
\qed


\begin {thebibliography} {99}

\bibitem{BR}{\sc C.~Benson and G.~Ratcliff}, On multiplicity free actions, Representations of real and $p$-adic groups, Lect. Notes Ser. Inst. Math. Sci. Natl. Univ. Singap., {\bf 2}, Singapore Univ. Press, Singapore, 2004, 221--304.


\bibitem{Bourbaki-Ch1-3}{\sc N.~Bourbaki}, Lie groups and Lie algebras, Chapters 1–-3, Translated from the French,  Reprint of the 1989 English translation, Springer-Verlag, Berlin, 1998. 


\bibitem{FH}{\sc W.~Fulton and J.~Harris}, Representation theory---A first course, Graduate Texts in Mathematics, {\bf 129}, Springer-Verlag, New York, 1991.  

\bibitem{Howe}{\sc R.~Howe}, Perspectives on invariant theory: Schur duality, multiplicity-free actions and beyond, 
The Schur lectures (1992) (Tel Aviv), Israel Math. Conf. Proc., {\bf 8}, Bar-Ilan Univ., 
Ramat Gan, 1995, 1--182. 


\bibitem{Humphreys}{\sc J. E.~Humphreys}, Introduction to Lie Algebras and Representation Theory, Graduate Texts in Mathematics, {\bf 9}, Springer-Verlag, New York-Berlin, 1972. 


\bibitem{Kac}{\sc V. G.~Kac}, Some remarks on nilpotent orbits, J. Algebra {\bf 64} (1980), 190--213.

\bibitem{Knop}{\sc F.~Knop}, Some remarks on multiplicity free spaces, Representation theories and algebraic geometry,  NATO Adv. Sci. Inst. Ser. C: Math. Phys. Sci., {\bf 514}, Kluwer Acad. Publ., Dordrecht, 1998, 301--317.  


\bibitem{Special1}{\sc K.~Nakamoto and Y.~Omoda}, 
Thick representations and dense representations I, 
Kodai Math. J. {\bf 42} (2019), 274--307.

\bibitem{panyushev}{\sc D. I.~Panyushev}, Properties of weight posets for weight
multiplicity free representations, Mosc. Math. J. {\bf 9} 
(2009), 867--883, 935--936. 

\end{thebibliography}
\end{document}